\documentclass[11pt]{article}
\usepackage{amsmath,amssymb,mathrsfs}
\usepackage{amsfonts}
\input xy
\xyoption{all}

\newcommand{\pf}{\noindent{\em Proof.}\ }
\newcommand{\qed}{\hfill $\Box$\\}

\renewcommand{\ni}{\noindent}

\newcommand{\ZZ}{\mathbb{Z}}

\newcommand{\tensor}{\otimes}

\newcommand{\dc}{\cap}
\newcommand{\rc}{\subset}

\newcommand{\isom}{\simeq}
\newcommand{\isoto}{\tilde{\to}}

\newcommand{\Spec}{\operatorname{Spec}}
\newcommand{\Spf}{\operatorname{Spf}}

\begin{document}

\title{Local structure of the moduli space of K3 surfaces
	over finite characteristic
	\footnote{This work is supported in part by Professor N. Yui's
		Discovery Grant from NSERC, Canada.}}
\author{\sc Jeng-Daw Yu}
\date{April 2007}
\maketitle

\begin{abstract}
Let $k$ be a perfect field of characteristic $p \geq 3$.
Let $X$ be a non-supersingular K3 surface over $k$,
and $\Psi$ the enlarged formal Brauer group associated to $X$.
We consider the deformation space of $X$.
In this note,
we show that the local moduli space
$\mathcal{M}^{\circ\circ}$ of $X$
with trivial associated deformation of $\Psi$
admits a natural $p$-divisible formal group structure.
\end{abstract}


\begin{center}
{\large 1. Introduction}
\end{center}

Throughout this note,
we let $p$ be an odd prime number.
Let $k$ be a perfect field of characteristic $p$.
Let $\sigma$ be the absolute frobenius automorphism on $k$.
Let $W$ be the ring of Witt vectors of $k$.
We also use $\sigma$ to denote the induced frobenius on $W$.\\

Fix a K3 surface $X$ over $k$ of finite height $h$
($1 \leq h \leq 10$).
Let $\mathfrak{Art}_k$ be the category of artinian local $k$-algebras.
Consider the formal deformation functor $\mathcal{M}$ of $X$
that sends every object $R$ in $\mathfrak{Art}_k$
to the isomorphism classes of formal deformations
$\mathfrak{X}$ of $X$ over $R$.
Then $\mathcal{M}$ is formally smooth of dimension 20 over $k$,
i.e. there is a non-canonical isomorphism
\[ \mathcal{M} \isom \Spf k[[t_1, \dots, t_{20}]]. \]
\\

If $X$ is ordinary, i.e. $h = 1$,
then it is known that
$\mathcal{M}$ admits a natural formal group structure.
More precisely,
let $\mathcal{D}_{<1}$ and $\mathcal{D}_{\leq 1}$
be the deformation space of the formal Brauer group $\Phi$
and the enlarged formal Brauer group $\Psi$
associated to $X$ respectively.
Then they are formally smooth of dimension
$h-1$ and $22-2h$ over $k$ respectively
and there are natural morphisms
\[ \xymatrix{ \mathcal{M} \ar[r]^{\alpha} & \mathcal{D}_{\leq 1}
	\ar[r]^{\beta} & \mathcal{D}_{< 1} } \]
which are formally smooth (\cite{NO}, Cor.(3.21)).
(The first map $\alpha$ sends a deformation of $X$ to its
associated enlarged formal Brauer group,
and the second map $\beta$ sends a deformation of $\Psi$
to its connected component.)
Then if $h = 1$,
the morphism $\alpha$ is an isomorphism (\cite{Nygaard}, Thm.1.6)
and
\[ \mathcal{D}_{\leq 1} \isom \Phi \tensor_{\ZZ_p} T^{\vee}_p(\Psi^{et}) \]
is a $p$-divisible formal group,
where $T^{\vee}_p(\Psi^{et})$ is the dual Tate module
of the \'etale quotient $\Psi^{et}$ of $\Psi$
(\cite{Chai-2slopes}, Prop.2.9).\\

Finally let $s = \Spec k$ and
$s \to \mathcal{D}_{< 1}$ and $s \to \mathcal{D}_{\leq 1}$
be the corresponding closed fibers.
Consider closed subspaces
$\mathcal{D}_{\leq 1}^{\circ}, \mathcal{M}^{\circ}$,
and $\mathcal{M}^{\circ\circ}$
defined in the following diagram with all squares cartesian
\[ \xymatrix{ {\mathcal{M}}^{\circ \circ} \ar[r] \ar[d] &
	{\mathcal{M}}^{\circ} \ar[r] \ar[d] & {\mathcal{M}} \ar[d] \\
	s \ar[r] & \mathcal{D}_{\leq 1}^{\circ} \ar[r] \ar[d] &
	\mathcal{D}_{\leq 1} \ar[d] \\
	& s \ar[r] & \mathcal{D}_{<1}. } \]
Thus, for example, for any $R$ in $\mathfrak{Art}_k$,
the set $\mathcal{M}^{\circ\circ}(R)$ consists of
isomorphism classes of deformations $\mathfrak{X}$ of $X$ over $R$
such that the associated enlarged formal Brauer group of $\mathfrak{X}$
is the trivial deformation of $\Psi$.\\

The aim of this note is to construct a natural formal group
structure on $\mathcal{M}^{\circ\circ}$
and show that this formal group is $p$-divisible.\\

In the case of local deformations of $p$-divisible groups
or abelian varieties over $k$,
Chai has discovered natural ``multi-extensions"
(called {\em cascades})
of $p$-divisible formal groups acting fully faithfully
on the completions of leaves,
which generalize the classical Serre-Tate coordinates
in the ordinary case
(\cite{Chai-survey}, Thm.7.7 and \cite{Chai-ICM}, \S 4.3;
see \cite{Chai-2slopes}, Thm.7.1 for two slopes case).
Our result can be regarded as a generalization
of Chai's results to the local deformations of K3 surfaces.
The main difference is that in our case,
there is no geometric object corresponding to the slope greater
than one part in the crystalline cohomology group
$H^2_{cris} (X/W)$ of degree two of the K3 surface $X$.
Instead of studying the whole deformation space
of the crystal $H^2_{cris} (X/W)$ and figuring out the subspace
coming from deformations of $X$ as in the approach of Chai,
we use the explicit description of crystals over $k[[t]]$,
and show that there is no $p$-torsion element
in $\mathcal{M}^{\circ\circ}(k[[t]])$.
Consequently $\mathcal{M}^{\circ\circ}$ is $p$-divisible.\\

It is a pleasure to express my appreciation to Prof. C.-L. Chai,
who shared the idea of his theory of canonical coordinates
and suggested the possible extension to the K3 surface case.
The author also thanks Prof. Yui for her encouragement and support.
\\

\begin{center}
{\large 2. Group structure on $\mathcal{M}^{\circ \circ}$}
\end{center}

Let $H = H^2_{cris} (X/W)$ be the crystalline cohomology
of the K3 surface $X$.
Consider the slope decomposition
\[ H = H_{<1} \oplus H_1 \oplus H_{>1}. \]
The slope less than one part $H_{<1}$
is the covariant Cartier module of the formal Brauer group $\Phi$ of $X$
and the slope one part $H_1$ corresponds to
the maximal \'etale quotient of $\Psi$.
The splitting is unique since $X$ is of finite height
and $k$ is perfect (\cite{Katz}, Thm.1.6.1).\\

In this section, a crystal $\mathcal{Q}$ over $R$ for an object
$k \to R$ in $\mathfrak{Art}_k$
will mean a crystal over $R$
relative to $W \to R$ through $k$
and the canonical divided power structure on $(W, pW)$.
\\

\ni{\bf 2.1.} The extension functor $\mathcal{E}$.\\

For $\iota: k \to R$ an artinian local $k$-algebra,
let
\[ \begin{array}{ccc}
\mathcal{H}_{<1} &=& \iota^* H_{<1} \\
\mathcal{H}_{>1} &=& \iota^* H_{>1}
\end{array} \]
be the {\em pull-back} crystals on $R$.
They are $F$-crystals over $R$.
Define $\mathcal{E}(R)$ to be
the set of isomorphism classes of $F$-crystals $\mathcal{P}$ on $R$
which is an extension of $\mathcal{H}_{>1}$ by $\mathcal{H}_{<1}$:
\[ 0 \to \mathcal{H}_{<1} \to \mathcal{P} \to \mathcal{H}_{>1} \to 0. \]
\\

For two objects $\mathcal{P}_1, \mathcal{P}_2$
in $\mathcal{E}(R)$,
define a binary operation
$\mathcal{P} = \mathcal{P}_1 + \mathcal{P}_2$
by the Baer sum as in the following diagram:

\[ \xymatrix{
0 \ar[r] &
\mathcal{H}_{<1} \oplus \mathcal{H}_{<1} \ar[r] &
\mathcal{P}_1 \oplus \mathcal{P}_2 \ar[r] &
\mathcal{H}_{>1} \oplus \mathcal{H}_{>1} \ar[r] & 0 \\
0 \ar[r] &
\mathcal{H}_{<1} \oplus \mathcal{H}_{<1}
	\ar@{=}[u] \ar[r] \ar[d]_{\Sigma} &
\mathcal{P}' \ar[u] \ar[r] \ar[d] &
\mathcal{H}_{>1} \ar[u]_{\Delta} \ar[r] \ar@{=}[d] & 0 \\
0 \ar[r] & \mathcal{H}_{<1} \ar[r] & \mathcal{P} \ar[r] &
	\mathcal{H}_{>1} \ar[r] & 0, } \]
where $\mathcal{P}'$ is the pull-back
of $\mathcal{P}_1 \oplus \mathcal{P}_2$
along the diagonal embedding $\Delta(x) = (x,x)$
and $\mathcal{P}$ is the push-out of $\mathcal{P}'$
along module sum $\Sigma(y,z) = y+z$
for sections $x,y,z$.
This defines the same $\mathcal{P}$ as one reverses
the order of pull-back and push-out:
\[ \xymatrix{
0 \ar[r] &
\mathcal{H}_{<1} \oplus \mathcal{H}_{<1} \ar[r] &
\mathcal{P}_1 \oplus \mathcal{P}_2 \ar[r] &
\mathcal{H}_{>1} \oplus \mathcal{H}_{>1} \ar[r] & 0 \\
0 \ar[r] &
\mathcal{H}_{<1}
	\ar@{=}[d] \ar[r] \ar@{<-}[u]^{\Sigma} &
\mathcal{P}'' \ar@{<-}[u] \ar[r] \ar@{<-}[d] &
\mathcal{H}_{>1} \oplus \mathcal{H}_{>1}
\ar@{<-}[d]^{\Delta} \ar[r] \ar@{=}[u] & 0 \\
0 \ar[r] & \mathcal{H}_{<1} \ar[r] & \mathcal{P} \ar[r] &
	\mathcal{H}_{>1} \ar[r] & 0. } \]
Here $\mathcal{P}''$ is the push-out along $\Sigma$
and $\mathcal{P}$ is the the pull-back along $\Delta$.\\

\ni{\bf Lemma.}
{\em Under the Baer sum binary operation,
the functor $\mathcal{E}$ is a group functor.}\\

\pf
This is a routine checking.
\qed

\ni{\bf 2.2.}
The map $\mathcal{M}^{\circ \circ} \to \mathcal{E}$.\\

If $\mathfrak{X}$ over $R$ is a deformation of $X$,
then the crystalline cohomology
\[ \mathcal{H} := H^2_{cris} \left( \mathfrak{X}/(R/W) \right) \]
is a deformation of $H$.
That is $\mathcal{H}$ is an $F$-crystal over $R$
such that the restriction $\mathcal{H}|_k$ to $k$
is the crystal $H$.
Let $\mathcal{H}_1 = \iota^* H_1$.\\

\ni{\bf Lemma.}
{\em Assume $\mathfrak{X} \in \mathcal{M}^{\circ \circ}(R)$,
Then there is a natural embedding
\[ \mathcal{H}_1 \to \mathcal{H} \]
and the cup product is perfect on this subspace $\mathcal{H}_1$.}\\

\pf
The triviality of the deformation
of the group $\Psi$ associated to $\mathfrak{X}$
(which is just $\Psi \times_k R$)
implies that there is an embedding of $F$-crystals over $R$
\[ \mathcal{H}_1 := \iota^* H_1 \to \mathcal{H} \]
by \cite{NO}, Thm.(3.20).
Together with \cite{NO}, Prop.(3.17),
one sees that the restriction of the cup-product on $\mathcal{H}$
to $\mathcal{H}_1$
is perfect.
\qed

\ni{\bf Lemma.}
{\em Let $\mathcal{P}$ be the orthogonal complement
of $\mathcal{H}_1$ with respect to the cup product pairing.
Then $\mathcal{P} \in \mathcal{E}(R)$.}\\

\pf
Still by \cite{NO}, Thm.(3.20),
the Dieudonn\'e module of the dual
of the associated formal Brauer group of $\mathfrak{X}$
gives the desired filtration.
\qed

Note the value of $\mathcal{P}$ at $R$ has an extra filtration
\[ {\rm Fil}_{\mathfrak{X}} \rc \mathcal{P}(R) \]
coming from the second Hodge filtration on the de Rham cohomology
$\mathcal{H}(R) = H^2_{dR} (\mathfrak{X}/R)$
of $\mathfrak{X}$.\\

\ni{\bf 2.3.}
Group law on $\mathcal{M}^{\circ \circ}$.\\

Let $(R, \mathfrak{m})$ be an object in $\mathfrak{Art}_k$,
where $\mathfrak{m}$ is the maximal ideal of $R$.
Let $Y, Z$ be deformations of $X$ in $\mathcal{M}^{\circ \circ}(R)$.
Assume $\mathfrak{m}^n = 0$ for some positive integer $n$.
Let $R_i = R/\mathfrak{m}^i$ and
$(Y_i, Z_i) = (Y, Z) \times_R R_i$.
As $(\mathfrak{m}^i)^2 = 0$ in $R_{i+1}$,
we could regard the natural quotient
$R_{i+1} \to R_i$ as a PD-thickening of $R_i$
with the trivial PD-structure on the kernel.\\

Let $\mathcal{P, Q}, \mathcal{P}_i, \mathcal{Q}_i$
be the corresponding deformations of the crystal $H$
associated to
$Y, Z, Y_i, Z_i$ respectively.
Define a new K3 surface $``Y + Z"$ over $R$
successively by defining K3 surfaces $(Y+Z)_i$ over $R_i$
in the following procedure:\\

(0) Given a K3 surface $\mathfrak{X}$ over $R_{i+1}$,
one gets a pair
$\left(\mathfrak{X} \times_{R_{i+1}} R_i, {\rm Fil}_{\mathfrak{X}} \rc
	H^2_{dR} (\mathfrak{X}/R_{i+1})\right)$
consisting of a K3 surface over $R_i$
and a rank one free and cofree $R_{i+1}$-submodule
${\rm Fil}_{\mathfrak{X}}$ of $H^2_{dR} (\mathfrak{X}/R_{i+1})$
coming from the Hodge filtration on the de Rham cohomology.
Note ${\rm Fil}_{\mathfrak{X}}$ is isotropic with respect to the cup product.

Under the above application,
we have a one-to-one correspondence
(see \cite{Deligne}, Th.2.1.11 and its proof)
\[ \{ \text{isomorphism classes of K3 surfaces over $R_{i+1}$} \} \]
\[ \updownarrow \]
\[ \left\{ \left(\mathfrak{X}/R_i, {\rm Fil} \rc
	H^2_{cris} \left(\mathfrak{X}/(R_i/W)\right)(R_{i+1}) \right) \right\}, \]
where the second set consists of
isomorphism classes of pairs
of a K3 surface $\mathfrak{X}$ over $R_i$
and an isotropic free and cofree submodule
of $H^2_{cris} \left(\mathfrak{X}/(R_i/W)\right)(R_{i+1})$
that lifts
${\rm Fil}_{\mathfrak{X}} \rc H^2_{dR} (\mathfrak{X}/R_i)$
via the restriction and the canonical identification
\[ H^2_{cris} \left(\mathfrak{X}/(R_i/W)\right)(R_{i+1}) \to
	H^2_{cris} \left(\mathfrak{X}/(R_i/W)\right)(R_i) =
	H^2_{dR} (\mathfrak{X}/R_i). \]
\\

(1) Over $R_1 = R/\mathfrak{m}$,
we have $Y_1 = Z_1 = X \times_k R/\mathfrak{m}$.
Define $(Y+Z)_1 = X \times_k R/\mathfrak{m}$.
Thus the new K3 surface $Y+Z$ will also be a deformation of $X$.\\

(2) Suppose we have defined $(Y+Z)_i$ over $R_i$.
Define $(Y+Z)_{i+1}$ over $R_{i+1}$ to be the K3 surface
corresponding to the pair
\[ \left( (Y+Z)_i, {\rm Fil} \right) \]
by the remark in step (0) above.
Here
\[ {\rm Fil} \rc (\mathcal{P}_i + \mathcal{Q}_i) (R_{i+1}) \rc
	[(\mathcal{P}_i + \mathcal{Q}_i) \oplus \mathcal{H}_1] (R_{i+1}) \]
is the induced filtration
from the filtrations
${\rm Fil}_{Y_{i+1}} \rc \mathcal{P}_i(R_{i+1})$ and
${\rm Fil}_{Z_{i+1}} \rc \mathcal{Q}_i(R_{i+1})$
associated to $Y_{i+1}$ and $Z_{i+1}$ respectively
under the Baer sum $\mathcal{P}_i + \mathcal{Q}_i$ in \S 2.1.\\

\ni{\bf Lemma.}
{\em Let $f: R \to A$ be a homomorphism in $\mathfrak{Art}_k$.
Then $f^*Y + f^*Z = f^*(Y + Z)$.}\\

\pf
This is trivial.
\qed


\ni{\bf Theorem.}
{\em The above binary operator defines
a natural group structure on the formal scheme $\mathcal{M}^{\circ\circ}$.
Furthermore this group $\mathcal{M}^{\circ\circ}$
is a formal group law of dimension $h-1$ over $k$.}\\

\pf
Since $\mathcal{M}^{\circ \circ} \isom k[[t_1, \dots, t_{h-1}]]$
is already pro-representable,
we only need to show that the group structure
on the tangent space
$\mathcal{M}^{\circ \circ}(k[\epsilon])$,
where $\epsilon^2 = 0$,
is the usual addition as in the vector space over $k$.

Since $k[\epsilon]$ is a PD-thichening of $k$,
under the natural identification as $k[\epsilon]$-modules
\[ H^2_{dR}(X/k) \tensor_k k[\epsilon] =
	H^2_{cris}(X/W) \tensor_W k[\epsilon]\ \isoto\
	H^2_{dR} (\mathfrak{X}/k[\epsilon]), \]
any lift $\mathfrak{X}$ of $X$ over $k[\epsilon]$
with trivial deformation of the enlarged formal Brauer group
is given by a map
\[ H^0 (X, \Omega^2) \to
	H^0 (X, \Omega^2)^{\perp} \dc \left(H_{<1} \tensor_W k\right) \]
(cf. \cite{Deligne}, proof of Th.2.1.11),
and the group structure reduces to the usual sum
(cf. see Section 3 below for explicit computation).
\qed
\\

\begin{center}
{\large 3. Divisibility}
\end{center}

In order to show that the formal group law on $\mathcal{M}^{\circ \circ}$
defined above is $p$-divisible,
it suffices to show that
there is no non-trivial $p$-torsion element in the Cartier module of
$\mathcal{M}^{\circ \circ}$.
For this purpose,
we may and do assume that $k$ is algebraically closed.
We will use the explicit describtion of crystals over $k[[t]]$
as in \cite{Katz}, \S 2.4.\\

Since $k$ is algebraically closed,
we can choose basis
$\{a_i\}$, $\{b_i\}$ for $H_{<1}$, $H_{>1}$
(regarded as free $W$-modules)
such that the absolute frobenius $F$ on them is given by
\[ F a_i = p^{(1-\delta_{hi})} a_{i+1} \]
\[ F b_i = p^{(1+\delta_{hi})} b_{i+1} \]
respectively and the cup product pairing $<\ ,\ >$ is given by
\[ < a_i, a_j > = 0 = < b_i, b_j > \]
\[ < a_i, b_j > = \delta_{ij} \]
for all $i, j$.
Here the index $i$ is to be understood as the residue modulo $h$.\\

\ni{\bf 3.1.} The trivial extension.\\

Let $\iota: k \to k[[t]]$ be the structure morphism.
As before, we write
$\mathcal{H}_{<1} = \iota^*H_{<1}$ and
$\mathcal{H}_{>1} = \iota^*H_{>1}$.
We fix a lift $\varphi: W[[t]] \to W[[t]]$
of the frobenius $\sigma: k[[t]] \to k[[t]].$
In the trivial deformation $\mathcal{H}_{<1} \oplus \mathcal{H}_{>1}$
of $H_{<1} \oplus H_{>1}$ to $k[[t]]$,
the (inverse images of the) elements $\{ a_i, b_i \}$ gives rise to a basis,
still denote them by $a_i, b_i$,
such that
\[ \nabla a_i = 0 = \nabla b_i \]
\[ F(\varphi)\varphi^* a_i = p^{(1-\delta_{hi})} a_{i+1} \]
\[ F(\varphi)\varphi^* b_i = p^{(1+\delta_{hi})} b_{i+1} \]
\[ < a_i, a_j > = 0 = < b_i, b_j > \]
\[ < a_i, b_j > = \delta_{ij}. \]
\\

If we change the basis which respects the filtration
\[ 0 \to \mathcal{H}_{<1} \to \mathcal{H}_{<1} \oplus \mathcal{H}_{>1}
	\to \mathcal{H}_{>1} \to 0 \]
and reduces to the direct sum
as tensoring $W$ over $W[[t]]$
by considering the new basis
$\{ a_i, b_i' \}$ of $\mathcal{H}_{<1} \oplus \mathcal{H}_{>1}$
where
\[ b_i' = b_i + \sum_{j=1}^h \alpha_{ij} a_j \]
with $\alpha_{ij} \in tW[[t]]$,
then
\[ \nabla b_i' = \sum_j d\alpha_{ij} \tensor a_j \]
and
\begin{eqnarray*}
F(\varphi)\varphi^* b_{i+1}' &=&
	p^{(1+\delta_{hi})} b_{i+1} +
	\sum_j \varphi^* \alpha_{ij} p^{(1-\delta_{hj})} a_{j+1} \\
	&=& p^{(1+\delta_{hi})} b_{i+1}' +
	\sum_j \left(p^{(1-\delta_{h, j-1})} \varphi^* \alpha_{i, j-1}
	- p^{(1+\delta_{hi})} \alpha_{i+1,j}\right) a_j.
\end{eqnarray*}
The cup product reads
\[ < a_i, a_j > = 0;\ < a_i, b_j' > = \delta_{ij};\
< b_i', b_j' > = \alpha_{ij} + \alpha_{ji} \]
for all $1 \leq i,j \leq h$.
These formulas give the conditions that a trivial extension
of $\mathcal{H}_{>1}$ by $\mathcal{H}_{<1}$
which deforms $H_{<1} \rc H_{<1} \oplus H_{>1}$
should satisfy.\\

\ni{\bf 3.2.} Freeness of $p$-torsion.\\

Now we are ready to show the following.\\

\ni{\bf Theorem.}
{\em The formal group law on $\mathcal{M}^{\circ\circ}$
defined in Section 2 is $p$-divisible.}\\

\pf
Suppose $\mathfrak{X} \in \mathcal{M}^{\circ\circ}(k[[t]])$
is a formal deformation of $X$ to $k[[t]]$ such that
the associated enlarged formal group of $\mathfrak{X}$
is the trivial deformation of $\Psi$.
As before we let $\mathcal{H}$ be the crystal over $k[[t]]$
given by the second crystalline cohomology of $\mathfrak{X}$
and $\mathcal{P}$ be the orthogonal complement of
$\iota^*H_1$ in $\mathcal{H}$.
Then one has an exact sequence
\[ 0 \to \mathcal{H}_{<1} \to \mathcal{P} \to
	\mathcal{H}_{>1} \to 0. \]
\\

Now suppose $\mathfrak{X}$ in $\mathcal{M}^{\circ \circ} (k[[t]])$
is $p$-torsion.
Choose elements $c_i$ in $\mathcal{P}$
that lift $b_i$ in $\mathcal{H}_{>1}$.
Write, for all $1 \leq i \leq h$,
\[ \nabla c_i = \sum_j \xi_{ij} \tensor a_j \]
\[ F(\varphi)\varphi^* c_i = p^{(1+\delta_{hi})} c_{i+1} +
	\sum_j v_{ij} a_j \]
($v_{ij} = v_{ij}(\varphi) \in tW[[t]]$) and
\[ < c_i, c_j > = (1 + \delta_{ij}) m_{ij} \]
\[ < a_i, c_j > = \delta_{ij} \]
for all $1 \leq i,j \leq h$.
Then by the discussion in \S 3.1,
there exist $\alpha_{ij} \in tW[[t]]$ such that
\begin{eqnarray}
p \xi_{ij} &=& d \alpha_{ij} \\
p v_{ij} &=&
	p^{(1-\delta_{h,j-1})} \varphi^* \alpha_{i,j-1}
	- p^{(1+\delta_{hi})} \alpha_{i+1,j} \\
p m_{ii} &=& \alpha_{ii} \\
p m_{ij} &=& \alpha_{ij} + \alpha_{ji} \hspace{15pt} i \neq j.
\end{eqnarray}
Since the crystal $\mathcal{H} = \mathcal{P} \oplus \mathcal{H}_1$
comes from a K3 surface $\mathfrak{X}$,
the frobenius $F$ on $\mathcal{H}$ is of rank 1 when modulo $p$
because $H^2(\mathfrak{X}, \mathcal{O})$ is of rank 1
(see \cite{Deligne}, (1.3.1.4)).
Thus
\begin{equation}
v_{ij} \equiv 0\
\text{(mod $p$) for all $1 \leq i \leq h$ and $j = 2, 3, \dots, h$}.
\end{equation}
Totally we have\\

$\bullet$
$\alpha_{ii} \equiv 0$ (mod $p$) by (3).\\

$\bullet$
If $\alpha_{ij} \equiv 0$, then
$\alpha_{ji} \equiv 0$ (mod $p$) by (4).\\

$\bullet$
For $j = 1, 2, \dots, h-1$, by (2),
\begin{eqnarray*}
pv_{i,j+1} &=& p^{(1-\delta_{h,j})} \varphi^* \alpha_{ij}
	- p^{(1+\delta_{hi})} \alpha_{i+1,j+1} \\
	&=& p\left(\varphi^* \alpha_{ij} -
	p^{\delta_{hi}} \alpha_{i+1,j+1}\right),
\end{eqnarray*}
i.e.
\begin{equation}
v_{i,j+1} = \varphi^* \alpha_{ij} -
	p^{\delta_{hi}} \alpha_{i+1,j+1}.
\end{equation}
Thus by (5), if $\alpha_{i+1,j+1} \equiv 0$,
then $\alpha_{ij} \equiv 0$
for $j = 1, 2, \dots, h-1$.\\

$\bullet$
Put $j = 1$ in (2):
\[ pv_{i1} = \varphi^* \alpha_{ih} -
	p^{(1+\delta_{hi})} \alpha_{i+1,1} \]
which implies
\[ \alpha_{ih} \equiv 0\ \text{(mod $p$)}. \]
Thus inductively via (6), 
we have $\alpha_{ij} \equiv 0$
for all $i,j$.
Therefore the deformation $\mathcal{P}$
(and hence $\mathcal{H}$) is trivial.\\

On the other hand,
if $\mathcal{P}$ is a trivial extension of
$\mathcal{H}_{>1}$ by $\mathcal{H}_{<1}$,
then the filtration on $\mathcal{H}(k[[t]])$
is the pull-back of the filtration on $H(k)$
since there is a unique rank 1 submodule on
$\mathcal{H}(k[[t]])$ such that
it is the
mod $p$ image of elements in $\mathcal{H}$
whose images under the frobenius $F$
are divisible by $p^2$.
Hence the deformation $\mathfrak{X}$ is trivial
by the remark in step (0) in \S 2.3.
Consequently $\mathcal{M}^{\circ \circ}$
is a $p$-divisible formal group.
\qed


\ni{\bf Proposition.}
{\em The formal $p$-divisible group $\mathcal{M}^{\circ\circ}$
is of frobenius slope $2/h$.}\\

\pf
This is a direct generalization of Prop.2.7 and Thm.2.8
in \cite{Chai-2slopes}.
One just needs to replace Barsotti-Tate groups over $k$
by crystals over $k$ in Prop.2.7.
As in Thm.2.8,
we consider the associated second crystalline cohomology
of the universal family
of K3 surfaces over $\mathcal{M}^{\circ\circ}$
instead of the universal extension $\mathcal{E}_{X,Y}$
over $\mathcal{DE}(X,Y)$ in \cite{Chai-2slopes}.
Since $H_{<1}$ and $H_{>1}$ have slopes
$1-1/h$ and $1+1/h$ respectively,
the difference $2/h$ gives the slope of
the $p$-divisible formal group $\mathcal{M}^{\circ\circ}$.
\qed
\\

\begin{center}
{\large 4. Remarks}
\end{center}

\ni{\bf 4.1.}
Over an algebraically closed field of characteristic $p$,
the only invariant of a one-dimestional $p$-divisible formal group
is its height
and similarly the only invariant for an \'etale $p$-divisible group
is its rank.
Thus for non-supersingular K3 surfaces,
the Newton polygon stratification coincides with
the leaf structure, i.e.
the foliation characterized by that in each leaf,
the associated crystalline cohomology groups
are geometrically constant.
Thus the results in this note confirm the expectation
that on each leaf in a good moduli space,
there would exist a fine structure from the extensions
of $p$-adic invariant.\\

\ni{\bf 4.2.}
In this aspect,
one should think of the subspace $\mathcal{D}_{<1}$ in \S 1
as the {\em normal direction} of a (non-supersingular)
Newton polygon stratum in the whole moduli space.\\

\ni{\bf 4.3.}
On the other hand,
the space $\mathcal{M}^{\circ}$ is a closed formal subscheme
in a {\em biextension} $\mathcal{D}^{\circ}$
sitting in the following fiber diagram
\[ \xymatrix{
\mathcal{M}^{\circ\circ} \ar[r] & \mathcal{D}^{\circ} \ar[d]^{\pi} \\
	& \mathcal{D}_{\leq 1}^{\circ} \times \mathcal{D}_{\geq 1}^{\circ}.} \]
Here $\mathcal{D}_{\geq 1}^{\circ}$ is the functor
that sends an object $R \in \mathfrak{Art}_k$
to isomorphism classes of $F$-crystals $\mathcal{S}$ over $R$
with
\[ \mathcal{S}_1 \rc \mathcal{S} \]
which deform the filtration
\[ H_1 \rc H_1 \oplus H_{>1}. \]
The functor $\mathcal{D}^{\circ}$ sends $R$
to isomorphism classes of $F$-crystals $\mathcal{T}$ over $R$
which deform
\[ H_{<1} \rc H_{<1} \oplus H_1 \rc H \]
such that the inverse image
$\pi^{-1}(s \times \mathcal{D}_{\geq 1}^{\circ})$
is the subgroup
$\mathcal{M}^{\circ\circ} \times \mathcal{D}_{\geq 1}^{\circ}$
of all possible deformations
$\mathcal{E} \times \mathcal{D}_{\geq 1}^{\circ}$.
By the duality of Dieudonn\'e modules and a Tate twist,
there is an isomorphism
$\iota: \mathcal{D}_{\leq 1}^{\circ} \to \mathcal{D}_{\geq 1}^{\circ}$.
Then the morphism
$\mathcal{M}^{\circ} \to \mathcal{D}_{\leq 1}^{\circ}$
is the pull-back of
$\mathcal{D}^{\circ} \to
	\mathcal{D}_{\leq 1}^{\circ} \times \mathcal{D}_{\geq 1}^{\circ}$
by
$id \times \iota: \mathcal{D}_{\leq 1}^{\circ} \to
	\mathcal{D}_{\leq 1}^{\circ} \times \mathcal{D}_{\geq 1}^{\circ}$.\\

\ni{\bf 4.4.}
The $p$-divisible formal group $\mathcal{M}^{\circ\circ}$
should only depend on the crystalline cohomology $H$ of $X$.
It would be interesting
to determine the Cartier/Dieudonn\'e module
of $\mathcal{M}^{\circ\circ}$
in terms of $H$.
For the case of $p$-divisible groups and abelian varieties,
see \cite{Chai-2slopes}, \S\S 4, 5, and 8.\\





\begin{thebibliography}{MM}

\bibitem[C1]{Chai-survey} C.-L. Chai,
	Hecke orbits on Siegel modular varieties.
	{\em Geometric methods in algebra and number theory},
	 Progr. Math., 235, pp. 71-107, 2005.

\bibitem[C2]{Chai-2slopes} C.-L. Chai,
	Canonical coordinates on leaves of $p$-divisible groups:
	The two-slope case (version 01/10/2005). Preprint.
	Available in {\tt http://www.math.upenn.edu/$\sim$chai/papers.html}.

\bibitem[C3]{Chai-ICM} C.-L. Chai,
	Hecke orbits as Shimura varieties in positive characteristic.
	Preprint.
	Available in {\tt http://www.math.upenn.edu/$\sim$chai/papers.html}.

\bibitem[D]{Deligne} P. Deligne,
		Cristaux ordinaires et coordonn\'ees canoniques.
		{\em Algebraic Surfaces (Orsay, 1976-78)}, LNM. 868,
		pp. 80-137, 1981.

\bibitem[K]{Katz} N. Katz, Slope filtration of $F$-crystals.
	{\em Journ\'ees de G\'eom\'etrie Alg\'ebrique de Rennes
	(Rennes, 1978), Vol. I}, Ast\'erisque 63, pp. 113-163, 1979.

\bibitem[N]{Nygaard} N. O. Nygaard,
	The Tate conjecture for ordinary K3 surfaces over finite fields.
	{\em Invent. Math.} 74 (1983), 213-237.

\bibitem[NO]{NO} N. O. Nygaard and A. Ogus,
		Tate's conjecture for K3 surfaces of finite height.
	{\em Ann. of Math.} 122 (1985), 461-507.

\end{thebibliography}
\end{document}